\documentclass[reqno,10pt]{amsart}
\usepackage{amssymb}
\usepackage[english,russian]{babel}
\usepackage[T2A]{fontenc}
\usepackage{color}

\usepackage{graphicx}
\usepackage{graphics}
\usepackage[colorlinks]{hyperref}
\usepackage{graphicx}

\usepackage{epstopdf}
\usepackage{epstopdf}
\usepackage[left=3cm,right=3cm,top=3cm,bottom=3cm]{geometry}

\numberwithin{equation}{section}

\newtheorem{theorem}{Theorem}[section]
\newtheorem{thm}[theorem]{Теорема}
\newtheorem{lem}[theorem]{Лемма}

\newtheorem{cor}[theorem]{Следствие}
\newtheorem{prop}[theorem]{Предложение}

\allowdisplaybreaks

\begin{document}

\title[]{О кривой критических показателей\\ нелинейных эллиптических уравнений \\в случае нулевой массы}

\author{Я.~Ш.~ Ильясов}\footnote{Работа выполнена при поддержки гранта РФФИ 14-01-00736-a}

\address{Я.Ш. Ильясов \newline
Институт математики, Уфимский научный центр, РАН\newline
Уфа, Чернышевского 112}
\email{ilyasov02@gmail.com}
\maketitle

\begin{abstract}
Рассматриваются полулинейные эллиптические уравнения в ограниченных и неограниченных областях вида: $-\Delta u =\lambda|u|^{p-2}u- |u|^{q-2}u$. В плоскости показателей нелинейностей $p\times q$ вводятся, так называемые, кривые критических показателей, отделяющих на этой плоскости  области с качественно отличительными свойствами рассматриваемых уравнений и соответствующих  параболических уравнений. Найдены новые условия разрешимости уравнений, устойчивости и не устойчивости стационарных решений,  существования глобальных решений параболических уравнений рассматриваемых во всем пространстве.

\end{abstract}


\section{Introduction}

В настоящей работе  рассматривается уравнение  
\begin{equation}
\label{s1} 
 -\Delta u =\lambda|u|^{p-2}u- \mu |u|^{q-2}u~~~ \mbox{in}~~ D.
 \end{equation}
Здесь  $p>0$, $q>0$, $\lambda \geq 0$, $\mu \geq 0$ $D$ -- одна из следующих областей: $D=\mathbb{R}^N$, $D=\Omega $, $D=\mathbb{R}^N\setminus \overline{\Omega}$, где $\Omega$ - звездная ограниченная область в $\mathbb{R}^N$ с $C^1$-гладкой границей $\partial \Omega$, $\mathrm{N}\geq 1$. Исследуются слабые решения $u \in W(D):=\mathcal{D}^{1,2}(D)\cap L^\infty(D)$ удовлетворяющие следующим граничным условиям
\par\noindent
при  $D=\Omega $, 
\begin{equation}\label{bD}
	u=0~~\mbox{на}~~\partial \Omega,
\end{equation}
при  $D=\mathbb{R}^N$,
\begin{equation}\label{bDI}
	|u(x)|\to 0~~\mbox{при}~~|x| \to \infty. 
\end{equation}
при $D=\mathbb{R}^N\setminus \overline{\Omega}$,
\begin{equation}\label{bI}
		u=0~~\mbox{на}~~\partial \Omega~~ \mbox{и}~~|u(x)|\to 0~~\mbox{при}~~|x| \to \infty.
\end{equation}
Краевые задачи \eqref{s1}- (\ref{bD}), \eqref{s1}- (\ref{bDI}), \eqref{s1}- (\ref{bI}) являются вариационными с  функционалом Эйлера-Лагранжа 
\begin{eqnarray}
E_\lambda(u)=\frac{1}{2}\int_{D} |\nabla u|^2\,dx -\lambda\frac{1}{p}\int_{D}
|u|^{p}\,dx +
 \frac{1}{q}\int_{D} |u|^{q}\,dx,  ~~u \in W(D).   \label{euler}
\end{eqnarray}
Следуя \cite{BerestL}, мы называем \eqref{s1}  уравнением соответствующим случаю нулевой массы. Такое уравнение представляет собой предельный случай семейства  уравнений с  ненулевыми массами, т.е. уравнений \eqref{s1}, в которых вместо $\Delta u$  рассматривается $(\Delta u-m u)$ при $m>0$.   
В случае $\mu=0$, уравнение \eqref{s1} с одним из граничных условий (\ref{bD}), (\ref{bDI}), (\ref{bI}) представляет собой пример  классической краевой задачи, обладающей критическим показателем нелинейности. Впервые, изучение нелинейных  задач с критическими показателями нелинейности,  были начаты  в 60-ые годы в работах Похожаева  и Фуджита \cite{Fujita, poh}. В работе \cite{Fujita}, Фуджита, рассматривая параболическую задачу $u_t=\Delta u +\lambda|u|^{p-2}u$ в $D=\mathbb{R}^N$, показал, что показатель нелинейности $p_F=\frac{2(N+1)}{N}$ является критическим, в том смысле, что при $p\in (1,p_F)$   параболическая задача  не имеет  неотрицательных глобальных решений, тогда как при $p>p_F$ такие решения возможны.  Похожаев в \cite{poh},  для  эллиптической задачи \eqref{s1}- (\ref{bDI}) 
при $\mu=0$ и $\mathrm{N}>2$, показал, что существование положительных решений  возможно только при  $p \in (2, 2^*]$, где $2^*$ -- критический показатель Соболева, ($2^*=\frac{2N}{N-2}$, при $\mathrm{N}\geq 3$, $2^*=+\infty$, при $\mathrm{N}=1,2$).  В настоящее время теория критических показателей нелинейностей является одной из центральных тем исследований в теории нелинейного анализа. Интерес к данной проблематики связан,  как с внутренней логикой развития теории  нелинейных дифференциальных уравнений, так и с востребованностью в математических результатах по этой теме  в прикладных исследованиях ( см. напр. \cite{bandle, Deng,Diaz-vol-1, galakV, galakMP, MPoh, Poh3} и приведенные там ссылки). 

В целом  критический показатель нелинейности, можно определить как значение показателя $p^*$,  которое задает  интервалы  такие, что  уравнения, рассматриваемые с показателями $p$ из этих интервалов, имеют качественно отличительное поведение. Однако, предоставляется, что задача о нахождении критического показателя нелинейности, должна быть изменена, если рассматривать \eqref{s1}  как семейство уравнений,  параметризованное  двумя показателями  $(p,q) \in \mathbb{R}^2$. В этом случае, возникает более общая  \textit{задача о нахождении кривой  критических показателей}, разделяющей плоскость  параметров $p\times q$ на области такими, что уравнение \eqref{s1} проявляет качественно отличительное поведение для каждой из областей, в которой находятся показатели  $(p,q)$. Исследования по данной задаче проводились в работе  \cite{DIH}, где  при $(p,q)\in \mathcal{E}_0:=\{(p,q) \in  \mathbb{R}^2:~ 1<q<p<2\}$ была найдена кривая  критических показателей нового типа $\mathcal{C}(N) \subset \mathcal{E}_0$, отделяющая  в $\mathcal{E}_0$ области, для показателей $(p,q)$ в каждой из которых ассоциированные  с \eqref{s1}-(\ref{bD}) начально-краевые параболические задачи могут иметь  только устойчивые или  неустойчивые основные состояния. 
В настоящей работе, мы развиваем исследования \cite{DIH}  применительно ко всему  квадранту $\mathcal{E}:=\{(p,q):~p>0,q>0\}$. 

Как и в \cite{DIH}, существенным в нашем подходе является использование функции  $\mathbb{R}^+ \ni r \mapsto E(ru)$ при $u \in W$, которую после работ Похожаева \cite{Poh1, Poh2} принято назвать \textit{функцией расслоений}.
Следует подчеркнуть, что  когда рассматриваются уравнения с зависимостью только от одного показателя нелинейности, например, при $\mu=0$ или  $\lambda=0$ в \eqref{s1}, то при каждом $u \in W\setminus 0$ функция расслоений  $E(ru)$ имеет единственную стационарную точку $r_u>0$: $dE(r_uu)/dr=0$, одного и того же типа при всех $u \in W\setminus 0$, т.е. либо $d^2E(r_uu)/dr^2<0$, либо $d^2E(r_uu)/dr^2>0$. Поведение функции расслоений $E(ru)$ существенно усложняется, если рассматривать \eqref{s1} в зависимости от двух показателей нелинейности $(p,q)$. В этом случае, $E(ru)$ может иметь две ненулевые стационарные точки  разных типов, или вообще их не иметь. Данная трудность, в нашем подходе, преодолевается тем, что для исследования разрешимости  задач, наряду с уравнением расслоений $dE(ru)/dr=0$ и тождеством Похожаева, как это делается при традиционном подходе (см. \cite{poh}), мы включаем в анализ также  уравнение $d^2E(ru)/dr^2=0$. Дополнительным преимуществом такого подхода является то, что это дает возможность  не только найти необходимые условия для существование  решений, но и определять к какому типу стационарной точки функции расслоений $E(ru)$ они могут соответствовать, т.е.  $d^2E(ru)/dr^2<0$, $d^2E(ru)/dr^2>0$ или $d^2E(ru)/dr^2=0$. В свою очередь, данные характеристики решений позволяют в дальнейшем исследовать их устойчивость для  соответствующих нестационарных задач  (см. ниже леммы \ref{lemUn}, \ref{lem4}, \ref{lem5}, \ref{lemF}). 

Центральную роль в данной работе играет следующее отображение 
\begin{equation}\label{Curv}
	d^*(p,q)=N (p-2)(q-2)-2pq, ~~(p,q) \in \mathbb{R}^2. 
\end{equation}
Множество $\mathcal{C}(p,q):=\{(p,q) \in \mathbb{R}^2:~d^*(p,q)=0\}$ мы называем \textit{кривой критических показателей}. Эта кривая, а также кривые 
$p=2$, $q=2$, $p=2^*$, $q=2^*$ и $p=q$ разделяют плоскость показателей нелинейности на области, как это показано на  Рис. 1, Рис. 2, Рис. 3. 
\begin{figure}[!h]
	\centering
	\includegraphics[width=0.5\linewidth]{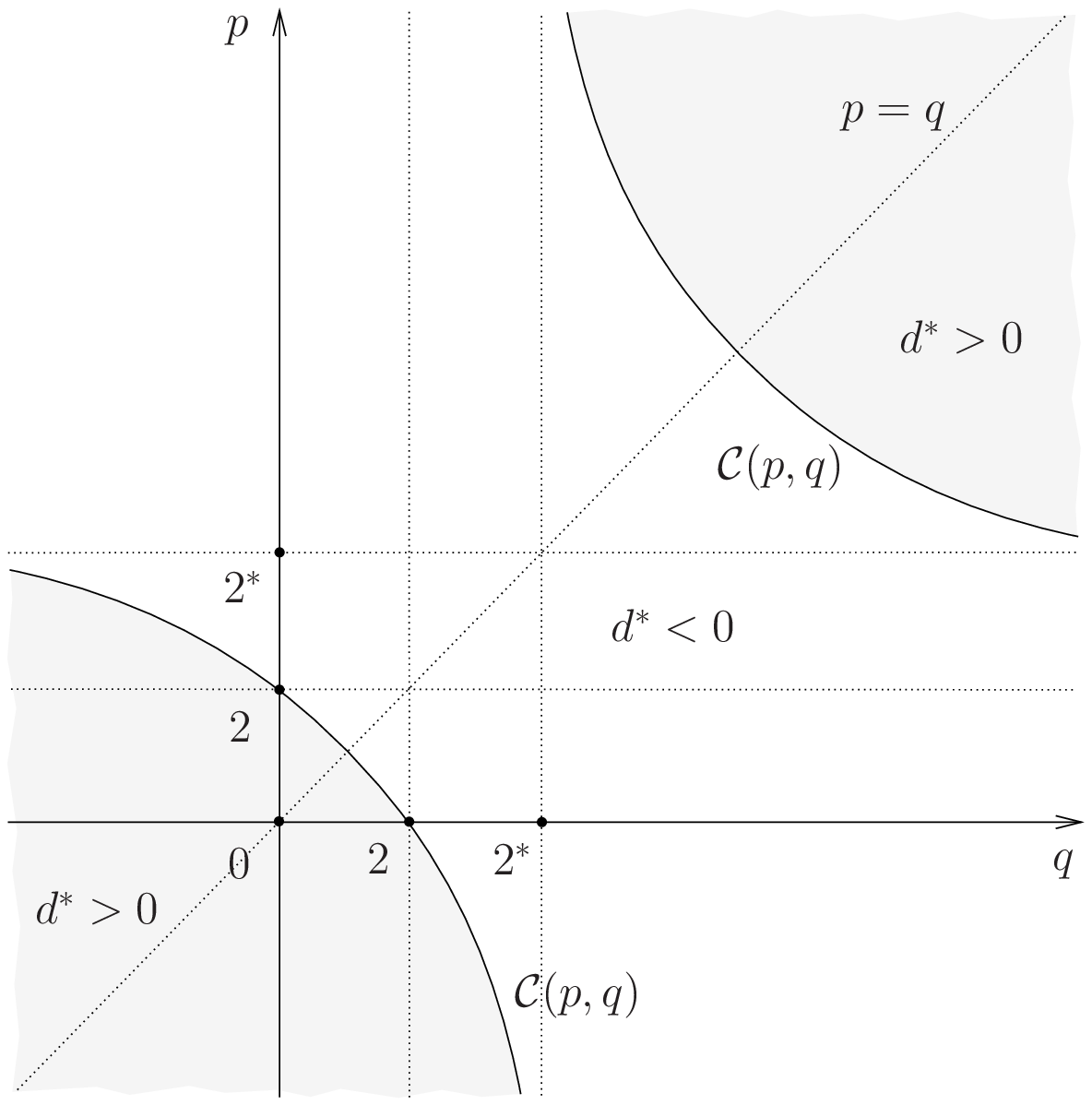}
	\caption{Кривая критических показателей при $N\geq 3$}
	\label{fig:Fig1}
\end{figure}
\begin{figure}[!ht]
\begin{minipage}[t]{0.45\linewidth}
\center{\includegraphics[scale=0.6]{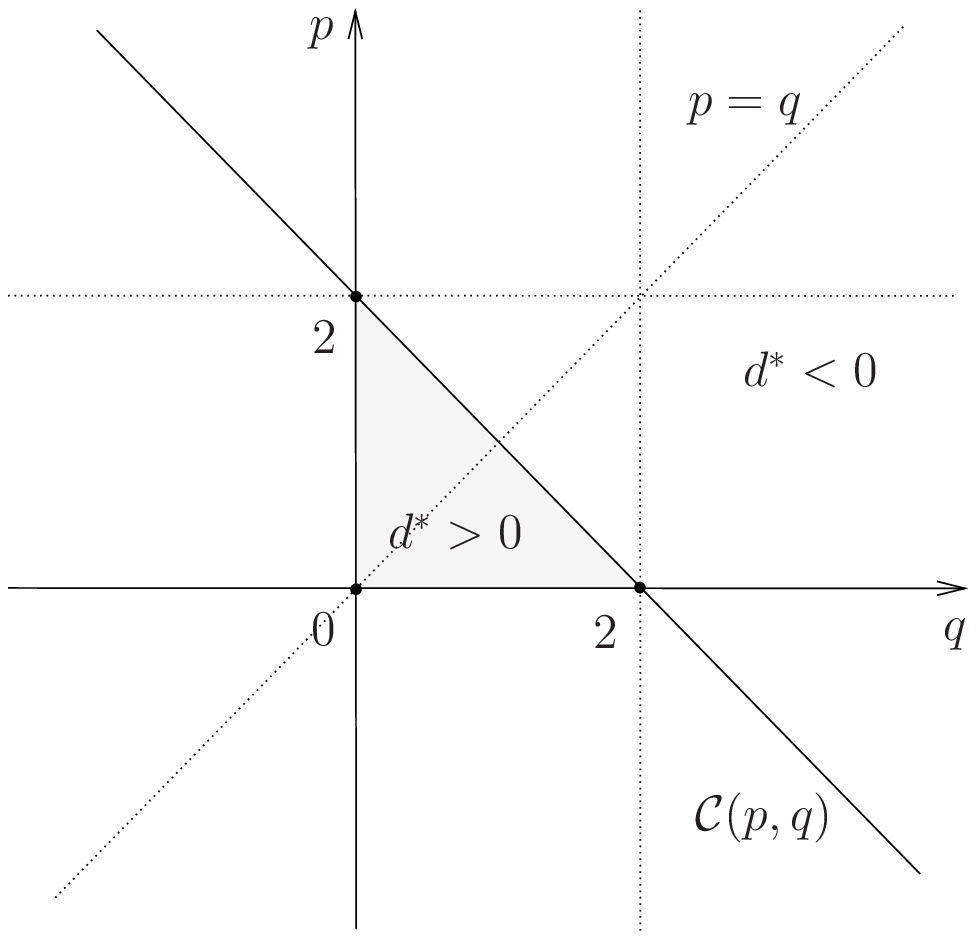}} \\
\caption{Кривая критических показателей при $N=2$}
\end{minipage}
\hfill
\begin{minipage}[t]{0.45\linewidth}
\center{\includegraphics[scale=0.6]{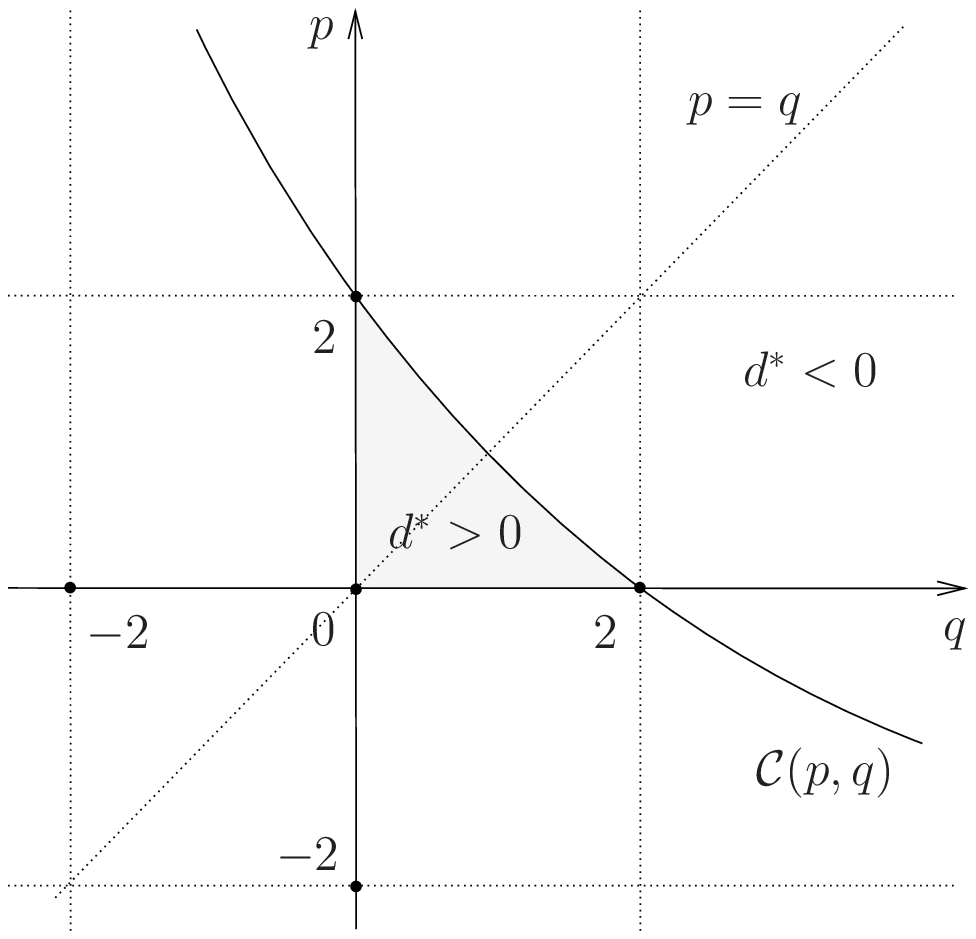}} \\
\caption{Кривая критических показателей при $N=1$}
\end{minipage}
\end{figure}

В основных результатах данной работы, исследуются свойства задач \eqref{s1} -(\ref{bD}), (\ref{bDI}), (\ref{bI}) в зависимости от того, какой из этих областей принадлежат показатели $(p,q)$. В разделе 2,  выводятся  необходимые условия существования решений задач \eqref{s1} -(\ref{bD}), (\ref{bDI}), (\ref{bI}) и, дается классификация решений, в зависимости от типа соответствующих им стационарных точек функции расслоений (см.  теоремы  \ref{Th1}, \ref{Th12}, \ref{Th3}). Как следствие этих результатов, мы находим ответ (в случае уравнений с нулевой массой) на проблему Штрауса (W.A. Strauss  \cite{Strauss}) о разрешимости задачи (\ref{s1})-(\ref{bDI})  при $2^*<q<p$ и $\mathrm{N}\geq 3$.
В разделе 3, мы исследуем существование решений. Основным  результатом здесь является теорема \ref{thmblO}. В разделе 4, используя качественные результаты о типах стационарных точек функции расслоений, мы получаем результаты о линейной неустойчивости стационарных решений параболических уравнений. В разделе 5, мы исследуем устойчивость решений параболических уравнений.   Здесь мы развиваем результаты \cite{DIH}, применительно к показателям нелинейности из квадранта $\mathcal{E}$. Отметим следующее, согласно теореме Деррика \cite{Derrick, Comech}, решения задачи (\ref{s1})-(\ref{bDI}), при $D=\mathbb{R}^N$, являются  линейно неустойчивыми стационарными состояниями для соответствующих параболических уравнений. Однако, этот результат справедлив только, если $p>2, q>2$, тогда как случаи $p,q \in (1,2]$ не подпадают под условия  теоремы Деррика. В данной работе, мы показываем, что при $p,q \in (1,2]$ справедлив результат, вообще говоря, противоположенный  утверждению теоремы Деррика.  В лемме \ref{lemF}, используя кривую критических показателей нелинейностей $\mathcal{C}(p,q)$, найдено подмножество $\mathcal{E}_u \subset (1,2]\times (1,2]$ обладающая следующими свойствами: при $(p,q) \in \mathcal{E}_u$, параболическая задача ассоциированная с  (\ref{s1})-(\ref{bDI})  обладает глобальными  решениями и устойчивыми (в заданном смысле) стационарными состояниями.

\section{Кривая критических показателей. Основные свойства.}

Уравнение \eqref{s1} достаточно исследовать в зависимости только от одного из параметров $\lambda>0$ или $\mu>0$. Действительно, делая замену переменных $\tilde{u}=tu$ в \eqref{s1}, например, с $t=(1/\mu)^{1/(q-2)}$, приходим к уравнению \eqref{s1} с $\tilde{\mu}=1$ и $\tilde{\lambda}=\lambda(1/\mu)^{(p-2)/(q-2)}$.  Более того, при $D=\mathbb{R}^N$, задача (\ref{s1})-(\ref{bDI}),   по существу, не зависит от обеих параметров $\mu, \lambda$. Действительно, если $u$ удовлетворяет (\ref{s1})-(\ref{bDI})  при некотором $\lambda>0$, то делая замену $v_{\tau,\sigma}=\tau u(\sigma x)$ c $\tau=\lambda^{1/(p-q)}$, $\sigma^2=\lambda^{(q-2)/(p-q)}$ получаем решение \eqref{s1} с $\tilde{\lambda}=1$.  
Таким образом, в дальнейшем, мы всюду полагаем  $\mu=1$. Кроме этого, рассматривая \eqref{s1} при $D=\mathbb{R}^N$, мы будем полагать, кроме специально оговоренных случаев, что и $\lambda=1$,  опуская в обозначениях индекс $\lambda$.

 Пусть  $D=\mathbb{R}^N$,  $D= \Omega$ или $D=\mathbb{R}^N\setminus \overline{\Omega}$, где $\Omega $ -- ограниченная область в $\mathbb{R}^N$ с $C^1$-гладкой границей $\partial \Omega$. Мы обозначаем $W(D):=\mathcal{D}^{1,2}(D)\cap L^\infty(D)$, где 
$\mathcal{D}^{1,2}(D)$ -- гильбертово пространство, задаваемое как пополнение $C^\infty_0(D)$ по норме $||w||_1=(\int_{D}|\nabla w|^2\,dx)^{1/2}$.  В этих обозначениях, краевые задачи (\ref{s1}) -(\ref{bD}), (\ref{bDI}), (\ref{bI})  записываются в единой форме
\begin{equation}\label{s1U}
\left\{ \begin{array}{clr}
	&	-\Delta u=\lambda|u|^{p-2}u-|u|^{q-2}u\\
	&~~u \in W(D),
		\end{array}\right.
\end{equation}
где равенство понимается в слабой форме.  Обозначим
\begin{equation*}
T(u):=\int_{D}|\nabla u|^{2}\,\mathrm{d}x,~~A(u)=\int_{D}|u|^p\,
\mathrm{d}x,~~B(u)=\int_{D}|u|^q\,
\mathrm{d}x, ~~u\in W(D).
\end{equation*}
Тогда $E_\lambda(u)=\frac{1}{2}T(u) -\lambda\frac{1}{p}A(u) +
\frac{1}{q}B(u)$.  Слабое решение $u \neq 0$ задачи \eqref{s1U} называется \textit{основным состоянием} (ground state \cite{BerestL}), если $E_\lambda(u)\leq E_\lambda(w)$ для любого другого слабого решения $w \in W(D)\setminus 0$ этой задачи. 
Рассмотрим  функцию расслоений
\begin{equation}\label{ffib}
	E_\lambda(ru)=\frac{r^2}{2}T(u) -\lambda\frac{r^p}{p}A(u) +
	\frac{r^q}{q}B(u), ~~r\in \mathbb{R}^+, ~~r>0,~~ u \in W(D)\setminus 0.
\end{equation}
Введем обозначения: $E'_\lambda(u):=\frac{d}{d r}E_\lambda(ru)|_{r=1}$, $E''_\lambda(u):=\frac{d^2}{d r^2}E_\lambda(ru)|_{r=1}$.  Отметим, что если  $u$ -- решение  задачи \eqref{s1U}, то $E'_\lambda(u)=0$, т.е. $r=1$  стационарная точка  функционала $E_\lambda(ru)$.  Легко видеть, что при каждом $u \in W(D)\setminus 0$ выполняется следующее: 
\par \noindent
(F1) \, если  $0<p<\min\{2,q\}$ или $\max\{2,q\}<p$, то при всех $\lambda\geq 0$, $E_\lambda(ru)$ имеет единственную ненулевую стационарную точку $r>0$, при этом $E''_\lambda(ru)>0$, если $1<p<\min\{2,q\}$, и   $E''_\lambda(ru)<0$, если $\max\{2,q\}<p$;
\par\noindent
(F2)\,  при $1<q<p<2$ или $2<p<q$, существует  такое $\lambda_u>0$, что при всех $\lambda\in  (0,\lambda_u)$, $E_\lambda(ru)$ не имеет ненулевых стационарных точек  $r>0$; при $\lambda=\lambda_u$ существует единственная ненулевая стационарная точка $r>0$, при этом $E''_\lambda(ru)=0$; при $\lambda\in (\lambda_u, +\infty)$ существуют две ненулевые стационарные  точки $r_{max}, r_{min}$ такие, что $E''_\lambda(r_{max}u)<0$, $E''_\lambda(r_{min}u)>0$.

Пусть $u\in W(D)$ слабое решение \eqref{s1U}.  Тогда по стандартной теорию регулярности решений эллиптических уравнений \cite{Lad} имееем   $u \in C^2(D)\cap C^{1,\kappa}(\overline{D})$ для $\kappa \in (0,1)$. Отсюда, $u$ удовлетворяет тождеству Похожаева \cite{poh} (см. также \cite{Takac_Ilyasov})
\begin{equation}
P_{\lambda }(u)+a\frac{1}{2N}\int_{\partial D }\left\vert \frac{%
\partial u}{\partial \nu }\right\vert ^{2}\,x\cdot \nu \,ds=0,  \label{poh}
\end{equation}%
где $\nu:=\nu(x)$ -- нормальный вектор к границе в точке $x\in D$, $a=0$, если $D=\mathbb{R}^N$, $a=1$, если $D=\Omega$ или $D=\mathbb{R}^N\setminus \overline{\Omega}$,
\begin{equation}
P_{\lambda }(u):=\frac{1}{2^*}T(u)-\lambda\frac{1}{%
p}A(u)+\frac{1}{q}B(u),~~~u\in W(D)
, \label{q landa}
\end{equation}
функция Похожаева. Здесь,  мы полагаем  $\frac{1}{2^*}=\frac{N-2}{2N}$, если  $\mathrm{N}\geq 3$, $\frac{1}{2^*}:=0$, если $\mathrm{N}=2$, $\frac{1}{2^*}=-\frac{1}{2}$, если  $\mathrm{N}=1$. Кроме этого, в дальнейшем мы будем использовать обозначения: $\frac{(q-2^*)}{2^*q}:=-\frac{1}{q}$, $\frac{(p-2^*)}{2^*p}:=-\frac{1}{p}$, если $\mathrm{N}=2$, и $\frac{(q-2^*)}{2^*q}:=-\frac{q+2}{q}$, $\frac{(p-2^*)}{2^*p}:=-\frac{p+2}{p}$, если  $\mathrm{N}=1$.
\begin{prop}\label{prop1}
 Пусть $u\in W(D)$ слабое решение \eqref{s1U}, $\Omega $ -- ограниченная область в $\mathbb{R}^N$, звездная относительно начала координат  $\mathbb{R}^{N}$, $\partial \Omega$ -- $C^1$-гладкая граница. Тогда
 \par\noindent
 (i) $P_{\lambda }(u)=0$, если $D=\mathbb{R}^N$;
  \par\noindent
  (ii) $P_{\lambda }(u)\leq 0$, если $D=\Omega$ и $u\geq 0$ в $D$. 
 \par\noindent
  (iii) $P_{\lambda }(u)\geq 0$, если $D=\mathbb{R}^N\setminus \overline{\Omega}$ и $u\geq 0$ в $D$.
  \end{prop}
	{\it Доказательство.} Утверждение (i) является непосредственным следствием формулы  (\ref{poh}). 
Отметим, что если $\Omega $ является звездной  областью  относительно начала координат  $\mathbb{R}^{N}$, то  $x\cdot \nu \geq 0$
при всех $x\in \partial \Omega $. Отсюда, поскольку для неотрицательных решений $\frac{
\partial u}{\partial \nu }\geq 0$ при всех $x\in \partial \Omega $, то справедливы (ii)-(iii).

\hspace*{\fill}\rule{3mm}{3mm}\\

Пусть $E_{\lambda }^{\prime }, E_{\lambda }^{\prime \prime }, P_{\lambda } \in \mathbb{R}$. Рассмотрим следующую систему уравнений  

\begin{equation}\label{sis1}
\left\{ 
\begin{array}{l}
\ T(u)-\lambda A(u)+ B(u)=E_{\lambda }^{\prime }, \\ 
 \\ 
\frac{1}{2^*}T(u)-\lambda\frac{1}{%
p}A(u)+ \frac{1}{q}B(u)=P_{\lambda },
\\
\\ 
\ T(u)-(p-1)\lambda A(u)+(q-1)B(u)=E_{\lambda }^{\prime \prime },
\end{array}\right.  
\end{equation}
где в качестве переменных рассматриваются $T(u), \lambda A(u), B(u) \in \mathbb{R}$. 
Детерминант этой системы равен
\begin{equation}
d=\frac{(q-p)(N(p-2)(q-2)-2pq)}{2Npq}=\frac{(q-p)}{2Npq}\cdot d^*(p,q),  \label{theta}
\end{equation}%
где $d^*$ задается по формуле (\ref{Curv}). 
При $d\neq 0$ и $E_{\lambda }^{\prime }=0$,  решение системы (\ref{sis1}) задается как:
\begin{equation}\label{Resis1}
\left\{
\begin{array}{l}
 T(u)=\frac{1}{d}\frac{(q-p)}{pq} E_{\lambda }^{\prime \prime }+ \frac{(q-p)}{d}P,\\ \\
\lambda A(u)=\frac{1}{d}\frac{(q-2^*)}{2^*q} E_{\lambda }^{\prime \prime }+ \frac{(q-2)}{d}P,\\ \\
		 B(u)=\frac{1}{d}\frac{(p-2^*)}{2^*p} E_{\lambda }^{\prime \prime }+ \frac{(p-2)}{d}P.
\end{array}\right.  
\end{equation}

В случае  $D=\mathbb{R}^N$, справедлив следующий результат о необходимых условиях 
\begin{thm}\label{Th1} Пусть $D=\mathbb{R}^N$,  $p\neq q$, $p>0, q>0$. Тогда

\par\noindent
$(1^o)$ для существования  ненулевого решения задачи (\ref{s1U}) необходимо   
	
	\begin{itemize}
		\item  если $\mathrm{N}\geq 3$, то $2^*<p<q$ или $0<q<p<2^*$;
		\item   если $\mathrm{N}=1,2$, то $0<q<p$.
	\end{itemize}
	
\par\noindent
$(2^o)$
если $u \in W(D)$ -- слабое ненулевое решение \eqref{s1U}, то  $E(u)>0$ и 
\begin{itemize}
\item $E^{\prime \prime }(u)=0~~\Leftrightarrow ~~d^*(p ,q )=0$,
	\item $E^{\prime \prime }(u)>0~~\Leftrightarrow ~~d^*(p ,q )>0$,
	\item $E^{\prime \prime }(u)<0~~\Leftrightarrow ~~d^*(p ,q )<0$;
	\end{itemize}
\par \noindent
\end{thm}
{\it Доказательство}.
 Пусть существует слабое ненулевое решение  $u$ задачи \eqref{s1U}. 
Тогда $E^{\prime }(u)=0$ и по предложению \ref{prop1}, $P(u)= 0$. Отсюда из (\ref{Resis1}) вытекает
\begin{equation}\label{Resis2}
\left\{ 
\begin{array}{l}
d \cdot T(u)=\frac{(q-p)}{pq} E^{\prime \prime }(u),\\ \\
\lambda d \cdot A(u)=\frac{(q-2^*)}{2^*q} E^{\prime \prime }(u),\\ \\
		 d\cdot B(u)=\frac{(p-2^*)}{2^*p} E^{\prime \prime }(u).
\end{array}\right.  
\end{equation}
Поскольку $T(u), A(u),B(u)>0$, то эти равенства возможны только, когда множители $(q-p)$, $(q-2^*)$, $(p-2^*)$ перед  $E^{\prime \prime }(u)$ имеют одинаковые знаки. Отсюда, при $\mathrm{N}\geq 3$, должно выполнятся $q>p$, $q>2^*$, $p>2^*$, либо $q<p$, $q<2^*$, $p<2^*$, а при $\mathrm{N}=1,2$, должно быть $q<p$. Отсюда, получаем $(1^o)$. 

Заметим, если $2^*<p<q$ и $\mathrm{N}\geq 3$, то $(q-p)>0$, $(q-2^*)>0$, $(p-2^*)>0$. Отсюда и из (\ref{Resis2}) вытекает, что знак $E^{\prime \prime }(u)$ совпадает со знаком $d^*(p,q)$. Легко видеть, что  аналогичное выполняется, если $0<q<p<2^*$ при $\mathrm{N}\geq 3$, или $0<q<p$ при  $\mathrm{N}=1,2$. Таким образом, получили $(2^o)$.  
\par\hspace*{\fill}\rule{3mm}{3mm}\\

В работе \cite{Strauss},  Штраус (W.A. Strauss)  сформулировал  вопрос о том, что может ли  задача (\ref{s1})--(\ref{bDI}) иметь решения при $2^*<q<p$ и $\mathrm{N}\geq 3$.  В случае нулевой массы, утверждение $(1^o)$ теоремы \ref{Th1} дает следующий ответ 
\begin{cor}
Задача (\ref{s1})--(\ref{bDI}), при $2^*<q<p$ и $\mathrm{N}\geq 3$, не имеет решений.
\end{cor}

В случае ограниченной области $D=\Omega$, справедлив следующий результат о необходимых условиях 
\begin{thm}\label{Th12} 
Пусть $p\neq q$,  $p>0, q>0$, $D=\Omega$, где $\Omega $ -- ограниченная область, звездная относительно начала координат  $\mathbb{R}^{N}$  с $C^1$-гладкой границей $\partial \Omega$. Тогда

\par\noindent
$(1^o)$  для существования  ненулевого неотрицательного решения задачи (\ref{s1U}) необходимо  
	\begin{itemize}
 \item  при $\mathrm{N}\geq 3:$ $0<p<q$, или  $0<q<p<2^*$; 
\end{itemize}

\par\noindent
 $(2^o)$	если $u \in W(D)$ -- слабое ненулевое неотрицательное решение \eqref{s1U}, то при $\mathrm{N}\geq 1$
\begin{itemize}
\item  $E^{\prime \prime }(u)>0$, если   $d^*(p,q)> 0$ или  $0<p<\min\{2,q\}$;
\item  $E^{\prime \prime }(u)<0$, если   $\max\{2,q\}<p<2^*$.
		\end{itemize}
	
\end{thm}
{\it Доказательство}.
 Предположим, что $u$ -- слабое ненулевое решение \eqref{s1U}. Пусть $\mathrm{N}\geq 3$. Тогда, т.к. $E_\lambda^{\prime }(u)=0$ и по предложению \ref{prop1}, $P_\lambda(u)\leq 0$, то из (\ref{Resis1}) вытекает
\begin{eqnarray}
&&d^*T(u)\leq\frac{1}{pq} E^{\prime \prime }(u), \label{Resis31}\\ 
&&\lambda\frac{d^*(q-p)}{(q-2)}A(u)\leq\frac{(q-2^*)}{2^*q(q-2)} E^{\prime \prime }(u), \label{Resis32}\\ 
		& &\frac{d^*(q-p)}{(p-2)} B(u)\leq\frac{(p-2^*)}{2^*p(p-2)} E^{\prime \prime }(u). \label{Resis33}
\end{eqnarray} 
Из (\ref{Resis31}) получаем, что если $d^*(p,q)>0$, то $E^{\prime \prime }(u)>0$. При этом неравенства (\ref{Resis32})-(\ref{Resis33}) не противоречивы, т.к. $d^*(p,q)>0$
влечет $\max\{p,q\}<2$ или $2^*<\min\{p,q\}$.

Рассмотрим  $\max\{2,q\}<p$.  Тогда, в силу (F1), $E^{\prime \prime }(u)<0$. Отсюда из (\ref{Resis31}) получаем $d^*(p,q)<0$. Тогда неравенство $d^*(p,q)(q-p)>0$ и   (\ref{Resis33}) влекут $p<2^*$. Таким образом, в этом случае  для существования решения необходимо  $\max\{2,q\}<p<2^*$. Легко видеть, что при $0<q<p<2$ неравенства (\ref{Resis31})-(\ref{Resis33}) непротиворечивы. Объединяя получаем, что, при $\mathrm{N}\geq 3$, в полуплоскости $0<q<p$ необходимым для существования решения \eqref{s1U} является условие $0<q<p<2^*$.

Рассмотрим случай $\mathrm{N}=2$. Тогда из (\ref{Resis1}) и, т.к. $P_\lambda(u)\leq 0$, имеем
\begin{equation}\label{Resis4}
\left\{ 
\begin{array}{l}
d^*T(u)\leq\frac{1}{pq} E^{\prime \prime }(u),\\ \\
\lambda\frac{d^*(q-p)}{(q-2)}A(u)\leq-\frac{1}{q(q-2)} E^{\prime \prime }(u),\\ \\
		\frac{d^*(q-p)}{(p-2)} B(u)\leq-\frac{1}{p(p-2)} E^{\prime \prime }(u).
\end{array}\right.  
\end{equation}
Из первого неравенства в (\ref{Resis4}) вытекает, что $E^{\prime \prime }(u)>0$, если $d^*(p,q)>0$.
При этом $d^*(p,q)>0$ влечет $p<2$ и $q<2$. Используя это, легко заключить, что эти неравенства непротиворечивы при всех рассматриваемых $0<p$ и $0<q$. При $\mathrm{N}=1$ рассуждения аналогичны. 
\par\hspace*{\fill}\rule{3mm}{3mm}\\

В случае  $D=\mathbb{R}^N\setminus \Omega$, справедлив следующий результат о необходимых условиях  
\begin{thm}\label{Th3}
Пусть $p\neq q$,  $p>0, q>0$, $D=\mathbb{R}^N\setminus \Omega$, где $\Omega$ -- ограниченная область  $\mathbb{R}^{N}$ с $C^1$-гладкой границей $\partial \Omega$, звездная  относительно начала координат. Тогда

\par\noindent
$(1^o)$  для существования  ненулевого неотрицательного решения задачи (\ref{s1U}) необходимо  
	\begin{itemize}
 \item  при $\mathrm{N}\geq 3$: $0<q<p$, или  $2^*<p<q$; 
\item  при $\mathrm{N}=1,2$: $0<q<p$, или  $2<p<q$;
\end{itemize}

\par\noindent
 $(2^o)$	если $u \in W(D)$ -- слабое ненулевое, неотрицательное решение \eqref{s1U}, то 
\begin{itemize}
\item при $\mathrm{N}\geq 3:$ $E^{\prime \prime }(u)<0$, если   $d^*(p,q)< 0$ и  $2^*<p<q$, или $d^*(p,q)< 0$ и $0<q<p$;
\item при $\mathrm{N}=1,2:$  $E^{\prime \prime }(u)<0$, если $2<p<q$, или $d^*(p,q)< 0$ и $0<q<p$.
		\end{itemize}

\end{thm}
{\it Доказательство}.
Предположим, что $u$ -- слабое ненулевое решение \eqref{s1U}. Тогда, при $\mathrm{N}\geq 3$, из (\ref{Resis1}), учитывая $P_\lambda(u)\geq 0$, имеем
\begin{eqnarray}
&&d^*T(u)\geq\frac{1}{pq} E^{\prime \prime }(u), \label{Resis31w}\\ 
&&\lambda\frac{d^*(q-p)}{(q-2)}A(u)\geq\frac{(q-2^*)}{2^*q(q-2)} E^{\prime \prime }(u), \label{Resis32w}\\ 
		& &\frac{d^*(q-p)}{(p-2)} B(u)\geq\frac{(p-2^*)}{2^*p(p-2)} E^{\prime \prime }(u). \label{Resis33w}
\end{eqnarray} 
Из (\ref{Resis31w}) получаем, что если $d^*(p,q)\leq 0$, то $E^{\prime \prime }(u)<0$. Рассмотрим $0<p<q$.
Поскольку, при $0<p<\min\{2,q\}$ мы имеем $E^{\prime \prime }(u)>0$, то если $d^*(p,q)\leq 0$ задача \eqref{s1U} не может иметь решений. Если же $2\leq p<q$, то (\ref{Resis32w}), (\ref{Resis33w}) возможны только при $2^*<p<q$. Легко видеть, что при $2<p<q$ и $d^*(p,q)>0$ неравенства (\ref{Resis31w})-(\ref{Resis33w}) допустимы, тогда как при $0<p<q<2$ и $d^*(p,q)>0$ они противоречивы. Таким образом, если $0<p<q$, то необходимо $2^*<p<q$ и если при этом $d^*(p,q)\leq 0$, то $E^{\prime \prime }(u)<0$. 
Анализ (\ref{Resis31w})-(\ref{Resis33w}), при  $0<q<p$, показывает, что эти неравенства  непротиворечивы.

Случаи $\mathrm{N}=1$ и $\mathrm{N}=2$ анализируются одинаково. Рассмотрим, как на примере, случай $\mathrm{N}=2$. Из (\ref{Resis1}), учитывая $P_\lambda(u)>0$, имеем
\begin{eqnarray}
&&d^*T(u)>\frac{1}{pq} E^{\prime \prime }(u), \label{Resis1w}\\ 
&&\lambda\frac{d^*(q-p)}{(q-2)}A(u)>-\frac{1}{q(q-2)} E^{\prime \prime }(u), \label{Resis2w}\\ 
		& &\frac{d^*(q-p)}{(p-2)} B(u)>-\frac{1}{p(p-2)} E^{\prime \prime }(u). \label{Resis3w}
\end{eqnarray} 
Анализ этих неравенств, при  $0<q<p$, показывает, что они непротиворечивы. 
Из (\ref{Resis1w}), если $d^*(p,q)\leq 0$, то $E^{\prime \prime }(u)<0$. 
Заметим, если $2<p<q$, то $d^*(p,q)\leq 0$ и следовательно левые части неравенств (\ref{Resis2w}), (\ref{Resis3w}) отрицательны, тогда как правые положительны. Получили противоречие. Легко видеть, что эти неравенства также противоречивы при $d^*(p,q)\leq 0$ и $0<p<\min\{2,q\}$. 
\par\hspace*{\fill}\rule{3mm}{3mm}\\

\section{Существование решений}
В этом параграфе мы исследуем существование решений задачи (\ref{s1U}) .

Существование решения (\ref{s1U}) при $D=\mathbb{R}^N$ вытекает из результатов полученных в \cite{BerestL, BerestLP, CortElgFelmer-2, Kaper1, Kaper2, Strauss}.  Мы просуммировали их в следующей
\begin{lem}\label{ThmEx1} Пусть $D=\mathbb{R}^N$. Тогда 
\begin{enumerate}
	\item если $\mathrm{N}=1,2$ и $1<q<p$, то существует  решение $u \in C^2(\mathbb{R}^N) $ задачи (\ref{s1U}) такое, что $u\geq 0$, $u(x)=u(r)$, $|x|=r$, $x \in \mathbb{R}^2$, $u'(r)<0$ при $r>0$. Кроме этого, если $2<p<q$, то  $u>0$ в $\mathbb{R}^N$.
	
	\item если $\mathrm{N}\geq 3$, $2^*<p<q$ или  $1<q<p<2^*$, то существует неотрицательное решение $u \in W(D)\cap C^2(\mathbb{R}^N)$ задачи (\ref{s1U}), являющейся основным состоянием. При этом, если $2^*<p<q$ или $2<q<p<2^*$, то $u>0$, функция $u$ сферически симметрична и монотонно убывает, т.е. $u(x)=u(r)$, $|x|=r$, $x \in \mathbb{R}^N$, $u'(r)<0$ при $r>0$;
	
	\end{enumerate}
\end{lem}
Доказательство утверждения (1) при $2<q<p$ вытекает, при $\mathrm{N}=2$, из теоремы 1.1 в \cite{BerestLP}, при $\mathrm{N}=1$,  из теоремы 5 в \cite{BerestL}. В случае $1<q<p\leq 2$, доказательство этого утверждения можно найти в \cite{DIH1, CortElgFelmer-2, Kaper1, Kaper2}.

Доказательство существования в утверждении (2)  вытекает из следующей  теоремы из \cite{BerestL} о достаточных условиях существования решений 
\begin{thm}(Berestycki-Lions)\label{thmbl}
Пусть $\mathrm{N}\geq 3$, $g:\mathbb{R} \to \mathbb{R}$ - непрерывная, нечетная функция такая, что
\begin{description}
	\item[(1)] $\overline{\lim}_{s\to +0}\frac{g(s)}{s^l}\leq 0$, где $l=2^*-1$;
	\item[(2)] $\exists\, \zeta>0$ такое, что $G(\zeta):=\int_0^\zeta g(s) ds>0$;
	\item[(3)] если $g(s)>0$ при всех $s>\zeta_0:=\inf\{\zeta>0:~G(\zeta)>0\}$, то 
	$$
	\lim_{s\to +\infty}\frac{g(s)}{s^l}=0.
	$$
\end{description}
Тогда существует неотрицательное решение $u \in \mathcal{D}^{1,2}(\mathbb{R}^N)$ задачи 
\begin{eqnarray}
	&&-\Delta u=g(u), \label{BL1}\\
	&& ~~u \in \mathcal{D}^{1,2}(\mathbb{R}^N). \label{BL2}
\end{eqnarray}
 При этом, $u \in C^2(\mathbb{R}^N)$. 
\end{thm}
Легко видеть, что при $1<\min\{p,q\}$, функция $g(s)=|s|^{p-2}s- |s|^{q-2}s$ непрерывна. Кроме этого, она  удовлетворяет условиям \textbf{(1)-(3)} теоремы \ref{thmbl} тогда и только тогда, когда $2^*<p<q$ или $1<q<p<2^*$. Отметим, что отсюда, в частности,  вытекает
\begin{cor} Пусть $g:\mathbb{R} \to \mathbb{R}$ - непрерывная, нечетная функция. Тогда условия \textbf{\rm (1)-(3)} теоремы  \ref{thmbl} являются необходимыми и достаточными   для существования решений задачи (\ref{BL1})-(\ref{BL2}).	
\end{cor}
Вторая часть утверждения (2) леммы \ref{ThmEx1} вытекает из принципа максимума для эллиптических уравнений
(см. напр. \cite{BerestL}).

Рассмотрим случай $D=\Omega$. В этом случае, мы будем рассматривать только показатели из множества $\max\{p,q\}<2^*$. Тогда по теореме Соболева справедливо вложение $\mathcal{D}^{1,2}(\Omega)\subset L^\gamma(D)$  для $\gamma \in (1,2^*)$. Поэтому, в качестве $W(D)$, можно рассматрвать $\mathcal{D}^{1,2}(\Omega)$. Мы будем строить решение (\ref{s1U}),  используя вариационную задачу с навязанным ограничением по методу многообразия Нехари:
\begin{equation}
\label{min1} 
\begin{array}{l}
\ E_\lambda(u) \to \min \\ 
~ u \in \mathcal{N}_\lambda.
%
\end{array}
\end{equation}
Здесь $\mathcal{N}_\lambda:=\{u \in W(D)\setminus 0: ~E'_\lambda(u)=0 \}$ называется многообразием Нехари. Обозначим 
\begin{equation}\label{minh} 
	\hat{E_\lambda}:=\min\{E_\lambda(u):~ u \in \mathcal{N}_\lambda\}.
\end{equation}
Отметим, поскольку любое слабое решение задачи (\ref{s1U}) принадлежит многообразию Нехари $\mathcal{N}_\lambda$, то  минимизирующая точка $u_\lambda \in \mathcal{N}_\lambda$ проблемы (\ref{min1}), удовлетворяющая уравнению (\ref{s1U}), является основным состоянием. 
 \begin{prop}\label{Lag2} Пусть $\lambda>0$ и $u_\lambda$ -- минимизирующая точка  задачи 
(\ref{min1}),  удовлетворяющая следующему условию
\begin{equation}\label{N}
E''_\lambda(u_\lambda)\neq  0.	
\end{equation}
Тогда $u_\lambda$ является слабым решением (\ref{s1U}), т.е. $D_uE_\lambda(u_\lambda)= 0$. 
 \end{prop}
{\it Доказательство}  см. напр. \cite{ Poh1, Poh2}.

Определим теперь значения $\lambda$, где $\mathcal{N}_\lambda \neq \emptyset$, а также выполняется условие (\ref{N}). Для этого мы воспользуемся методом нелинейно-обобщенных отношений Релея \cite{ilyaReil}. 
Пусть  $u \in \mathcal{D}^{1,2}(\Omega) \setminus 0$, $1<\min\{p,q\}$, $\max\{p,q\}<2^*$, тогда корректно определено следующее отношения Релея 
$$
	R(u)=\frac{\int_{D} |\nabla u|^{2} dx+\int_{D} |u|^{q} dx}{\int_{D} |u|^{p} dx}.
$$
Рассмотрим  
\begin{equation}\label{RaylCC}
	R(r u)=\frac{r^{2-p}\int_{D} |\nabla u|^{2} dx+r^{q-p}\int_{D} |u|^{q} dx}{\int_{D} |u|^{p} dx},~~u \in \mathcal{D}^{1,2}(\Omega) \setminus 0,  ~r>0.
\end{equation}
Дифференцируя по $r$ эту функцию, получаем
\begin{equation}
	\frac{\partial}{\partial r} R(r u)=\frac{(2-p)r^{1-p}\int_{D} |\nabla u|^{2} dx+(q-p)r^{q-p-1}\int_{D} |u|^{q} dx}{\int_{D}| u|^{p} dx}.
\end{equation}
Легко видеть, что    уравнение  $\frac{\partial}{\partial r} R(r u)=0$ имеет решение только при $1<q<p<2$ и $2<p<q$, при этом, решение единственно и задается по формуле 
\begin{equation}\label{tmaxcc}
	r_{min}(u)=\left(\frac{(p-2)\int_{D}|\nabla u|^{2} dx}{(q-p)\int_{D} |u|^{q} dx}\right)^\frac{1}{q-2}.  
\end{equation}
Подставляя $r_{min}(u)$ в  $R(r u)$, получаем следующее нелинейно-обобщенное отношений Релея 
\begin{equation}\label{LambCc}
	\lambda(u):=R(r_{min}(u) u)=c_{p,q}\cdot
	 \frac{(\int_{D}|\nabla u|^{2} dx)^{\frac{q-p}{q-2}}(\int_{D} |u|^{q} dx)^{\frac{p-2}{q-2}}}{\int_{D} |u|^{p} dx } 
\end{equation}
где
$$
c_{p,q}=\frac{(p+q-4)}{(q-2)} \left( \frac{(q-2)}{(p-2)}
	\right)^{\frac{p-2}{q-2}} 
$$
Отсюда получаем следующие критическое значение \cite{ilyaReil}:
\begin{equation}\label{lambz}
	\lambda_{(p,q)}=c_{p,q}\cdot \inf_{u \in \mathcal{D}^{1,2}(D) \setminus 0} \frac{(\int_{D}|\nabla u|^{2} dx)^{\frac{q-p}{q-2}}(\int_{D} |u|^{q} dx)^{\frac{p-2}{q-2}}}{\int_{D} |u|^{p} dx }. 
\end{equation}
Используя те же рассуждения, для следующего отношения Релея 
$$
	R_E(u)=\frac{\frac{1}{2}\int_{D} |\nabla u|^{2} dx+\frac{1}{q}\int_{D} |u|^{q} dx}{\frac{1}{p}\int_{D} |u|^{p} dx},
$$
получаем  
\begin{equation}\label{LambCc}
	\lambda_E(u):=R_E(r_{min}(u) u)=c_{p,q}'c_{p,q}\cdot
	 \frac{(\int_{D}|\nabla u|^{2} dx)^{\frac{q-p}{q-2}}(\int_{D} |u|^{q} dx)^{\frac{p-2}{q-2}}}{\int_{D} |u|^{p} dx } 
\end{equation}
где $c_{p,q}'=(p/2) \left( 2/q
	\right)^{\frac{p-2}{q-2}}$.
Отсюда мы можем ввести: $\lambda_{E,(p,q)}=\inf_{u \in \mathcal{D}^{1,2}(D) \setminus 0} \lambda_E(u)$. Очевидно 
$\lambda_{E,(p,q)}=c_{p,q}'\lambda_{(p,q)}$ и, при этом, 
\begin{lem} \label{Rl} Пусть $1<q<p<2$ или $2<p<q<2^*$, $D=\Omega $, где $\Omega$ -- ограниченная область в $\mathbb{R}^{N}$ с $C^1$-гладкой границей $\partial \Omega$. Тогда
$0<\lambda_{(p,q)}<+\infty$,  $\lambda_{(p,q)}<\lambda_{E,(p,q)}$. При этом, $\mathcal{N}_\lambda \neq \emptyset$ тогда и только тогда, когда $\lambda>\lambda_{(p,q)}$ и, если $\lambda>\lambda_{E,(p,q)}$, то найдется  $u \in \mathcal{D}^{1,2}(D) \setminus 0$	такое, что $E_\lambda(u)<0$.
\end{lem}
{\it Доказательство } Оценка $0<\lambda_{(p,q)}$ легко выводится, используя неравенства Гельдера и Соболева. Неравенство $\lambda_{(p,q)}<\lambda_{E,(p,q)}$ следует из того, что $c_{p,q}'>1$. 
Для доказательство оставшейся части, достаточно заметить, что  $\mathcal{N}_\lambda=\{u \in W(D)\setminus 0: ~R(u)=\lambda \}$, $\forall \lambda>0$ и, что если $R_E(u)<\lambda$ для некоторого $u \in \mathcal{D}^{1,2}(\Omega) \setminus 0$, то $E_\lambda(u)<0$. 

\hspace*{\fill}\rule{3mm}{3mm}\\

Докажем следующий основной результат о существовании решений                         

\begin{thm}\label{thmblO}
Пусть $N\geq 1$, $D=\Omega $, где $\Omega$ -- ограниченная область в $\mathbb{R}^{N}$ с $C^1$-гладкой границей $\partial \Omega$. Тогда

\begin{enumerate}
	\item если $1<p<\min\{2,q\}$ или $1<q$ и $\max\{2,q\}<p<2^*$,  то при всех $\lambda\in (0, +\infty)$ существует решение $u \in W(D)$ задачи (\ref{s1U}) такое, что $u$ -- основное состояние, $u \in C^2(D)\cap C^{1,\kappa}(\overline{D})$ для некоторого $\kappa \in (0,1)$ и  $u_\lambda>0$ в $D$. 
	\item  если $1<q<p\leq 2$ или $2<p<q<2^*$,  то задача (\ref{s1U}) не имеет решений при $0<\lambda <\lambda_{(p,q)}$. Если 	$\lambda >\lambda_{E,(p,q)}$, то  существует решение $u_\lambda \in W(D)$ задачи (\ref{s1U}). При этом, $u_\lambda$ является основным состоянием, $u_\lambda \in C^2(D)\cap C^{1,\kappa}(\overline{D})$ для некоторого $\kappa \in (0,1)$, $u_\lambda \geq 0$ в $D$ и $E_\lambda(u_\lambda)<0$, $E''_\lambda(u_\lambda)>0$, $\forall \lambda >\lambda_{E,(p,q)}$. Кроме этого, если $2<p<q<2^*$, то  $u_\lambda>0$ в $D$.
	\end{enumerate}
 
\end{thm}
\textit{Доказательство.} Отметим,  поскольку $E_\lambda(u)=E_\lambda(|u|)$ и $|u|\in \mathcal{N}_\lambda$, $\forall u \in \mathcal{N}_\lambda$, то из существования  минимизирующей точки $u_\lambda$ задачи 
(\ref{min1}) вытекает, что $|u_\lambda|$ является также  минимизирующей. При этом, условие (\ref{N}) очевидно также сохраняется. Таким образом, доказывая существование минимизирующей точка $u_\lambda$ задачи (\ref{min1}) удовлетворяющей (\ref{N}), мы получаем существование слабого неотрицательного решения (\ref{s1U}). Далее, если $u_\lambda$ слабое решение задачи (\ref{min1}), то теория регулярности решений эллиптических уравнений (см.
\cite{Lad}) влечет, что   $u_\lambda \in C^2(D)\cap C^{1,\kappa}(\overline{D})$  для некоторого $\kappa \in (0,1)$. Кроме этого, при  $1<p<\min\{2,q\}$, $1<\max\{2,q\}<p<2^*$ или $2<p<q<2^*$ из принципа максимума для эллиптических краевых задач вытекает \cite{giltrud, Vask}, что $u_\lambda>0$ в $D$. Таким образом, для доказательства 
утверждений (1)-(2), Теоремы \ref{thmbl} достаточно найти (при соответствующих $\lambda$) минимизирующую точку  задачи (\ref{min1}), удовлетворяющую условию (\ref{N}).

\textsl{Доказательство утверждения} (1). Заметим,  если  $1<p<\min\{2,q\}$, то  $E''_\lambda(u)>0$, $\forall u \in \mathcal{N}_\lambda$ и, если $1<q$ и $\max\{2,q\}<p<2^*$, то $E''_\lambda(u)<0$, $\forall u \in \mathcal{N}_\lambda$. Таким образом, в этих случаях, условие (\ref{N}) всегда выполняется и для доказательства утверждения (1) достаточно показать, что существует минимизирующая точка  задачи 
(\ref{min1}).

Пусть $u_m \in \mathcal{N}_\lambda$, $m=1,2,...$ -- минимизирующая последовательность проблемы (\ref{min1}), т. е. $E_\lambda(u_m) \to \hat{E_\lambda}$ при $m \to \infty$. 

Рассмотрим сначала случай $1<p<\min\{2,q\}$. Используя неравенства Гельдера и Соболева, выводится 
$$
E_\lambda(u)\geq \max\{ \frac{1}{2}T(u)- C\frac{1}{p}T(u)^{p/2} , - C\frac{1}{p}B(u)^{p/q} +
\frac{1}{q}B(u)\}
$$
где $0<C<+\infty$ не зависит от $u \in  W(D)$. Отсюда вытекает, что 
$E_\lambda(u)$ является коэрцитивным функционалом на $W(D)$, и следовательно найдется подпоследовательность, снова обозначаемой $(u_m )$ такая, что   $u_m \rightharpoondown u_\lambda \in W(D)$ слабо в  $W(D)$  и по теореме Соболева $u_m \to u_\lambda$ сильно в  $L^\gamma(D)$  для $\gamma \in (1,2^*)$. Отсюда, т.к. $p \in (1,2)$ вытекает, что $A(u_m) \to A(u_\lambda)$ и следовательно
\begin{align}
	&E_\lambda(u_\lambda) \leq \liminf_{m \to \infty}E_\lambda(u_m)=\hat{E}_\lambda \label{sh1}\\
	&E'_\lambda(u_\lambda) \leq \liminf_{m \to \infty}E'_\lambda(u_m)=0.\label{sh2}
\end{align}
Поскольку $E_\lambda(u_m) <0$, $m=1,2,...$, то $\hat{E}_\lambda<0$, и следовательно $u_\lambda \neq 0$. Легко видеть, что если в одном из неравенств (\ref{sh1})-(\ref{sh2}) выполняется равенство, то $u_\lambda$ является минимизирующей точкой  задачи 
(\ref{min1}).  Предположим, что $E'_\lambda(u_\lambda)<0$. Тогда найдется такое $r>1$, что 
$E'_\lambda(ru_\lambda)=0$ и $E_\lambda(ru_\lambda)<\hat{E}_\lambda$, что является противоречием. Таким образом, получаем $E_\lambda(u_\lambda)=\hat{E}_\lambda$ и $E'_\lambda(u_\lambda)=0$, т.е.  $u_\lambda$ является минимизирующей точкой  задачи 
(\ref{min1}). 

Рассмотрим теперь случай, когда $1<q$ и $\max\{2,q\}<p<2^*$. В этом случае, $E_\lambda(u)$ является коэрцитивным функционалом на $\mathcal{N}_\lambda$. Действительно, если  $u\in \mathcal{N}_\lambda$, то $E_\lambda(u)=\frac{p-2}{2p}T(u)+\frac{p-q}{pq}B(u) \to \infty$ при $\|u\|_{W} \to \infty$.  Отсюда, как и выше, выводится, что существует предельная функция $u_\lambda \in W(D)$, для которой выполняются (\ref{sh1})-(\ref{sh2}). Покажем, что $u_\lambda \neq 0$. Предположим противное. Тогда $r_m:=||u_m||_1 \to 0$ и
$$
0=T(v_m)-\lambda r_m^{p-2}A(v_m)+ r_m^{q-2}B(v_m)\to 1~~~\mbox{при}~~~m \to \infty,
$$
где $v_m=u_m/||u_m||_1$, $m=1,...$ . Получили противоречие. Покажем, что в (\ref{sh1})-(\ref{sh2}) справедливы равенства. Предположим противное. Тогда $E'_\lambda(u_\lambda)<0$ и найдется $r_0\in (0,1)$ такое, что $E'_\lambda(r_0u_\lambda)=0$. Поскольку $r_0<1$, а $r=1$ точка максимума функции $E_\lambda(ru_m)$, $m=1,...$, то 
$$
E_\lambda(r_0u_\lambda)\leq \liminf_{m \to \infty}E_\lambda(r_0u_m)\leq \liminf_{m \to \infty}E_\lambda(u_m)=\hat{E}_\lambda 
$$ 
Отсюда учитывая, что $r_0u_\lambda \in \mathcal{N}_\lambda$, получаем противоречие. Утверждение (1) доказано.


 Утверждение (2) при $1<q<p\leq 2$ доказано в \cite{DIH} (см. также \cite{CortElgFelmer-2, DIH1,Kaper1}). Докажем  это утверждение  при  $2<p<q<2^*$. Отсутствие решений задачи (\ref{s1U})  
при $0<\lambda <\lambda_{(p,q)}$ вытекает из  леммы \ref{Rl}. Действительно любое решение принадлежит многообразию Нехари, а по лемме \ref{Rl} мы имеем $\mathcal{N}_\lambda=\emptyset$, если $0<\lambda <\lambda_{(p,q)}$.

\begin{lem}\label{le1e} Пусть $2<p< q<2^*$ и $\lambda>\lambda_{(p,q)}$. Тогда существует  минимизирующая точка $u_\lambda \in \mathcal{N}_\lambda$ проблемы (\ref{min1}). 
\end{lem}
\textit{Доказательство.}
Пусть $\lambda>\lambda_{(p,q)}$. Тогда, по лемме \ref{Rl},  $\mathcal{N}_\lambda$ не пусто.
Пусть $u_m \in \mathcal{N}_\lambda$, $m=1,2,...$ -- минимизирующая последовательность проблемы (\ref{min1}), т. е. $E_\lambda(u_m) \to \hat{E_\lambda}$ при $m \to \infty$.  
Покажем, что $u_m$
ограниченно в $\mathcal{D}^{1,2}(\Omega)$.   Используя неравенство Гельдера, имеем
$$
E_\lambda(u_m)\geq -\lambda\frac{1}{p} A(u_m)	+\frac{1}{q}c_{\Omega}\,(A(u_m))^{q/p}, ~~m=1,...,
$$
где $0<c_{\Omega}<+\infty$ не зависит от $m$. Отсюда и, поскольку $q>p$, вытекает, что $A(u_m)$ ограниченно, что, в свою очередь, в силу $E'_\lambda(u_m)=0$, влечет ограниченность $T(u_m)$. Получили требуемое. Кроме этого, мы выяснили, что $\hat{E_\lambda}>-\infty$.  Рассмотрим $v_m=u_m/||u_m||_1$, $r_m:=||u_m||_1$, $m=1,2,...$ . Поскольку $r_m$ ограниченно, то, без ограничения общности, можем считать, что $r_m \to r_0$ при $m\to \infty$ для некоторого $r_0\geq 0$. Заметим, при $m=1,...$, справедливо 
	\begin{equation}\label{EQ}
		\frac{q-2}{2q}T(v_m)-\lambda r_m^{p-2}\frac{q-p}{pq}A(v_m) = \frac{E_\lambda(v_m)}{r_m^2}.
	\end{equation}
Отсюда, учитывая $-\infty<\hat{E_\lambda}$, получаем $r_0\neq 0$. 
 Поскольку $||v_m||=1$, $m=1,2,...$, то по теореме Эберлейна-Шмулияна
и теореме  вложения Соболева найдется подпоследовательность, снова обозначаемая $(v_m)$, такая, что $v_m \rightharpoondown v_0$ слабо в $\mathcal{D}^{1,2}(\Omega)$ и сильно $v_m \to v_0$ в $L_p(\Omega)$ и $L_q(\Omega)$  для некоторого $v_0\in \mathcal{D}^{1,2}(\Omega)$. Отсюда, рассуждая от противного с  использованием (\ref{EQ}),  выводится, что $v_0 \neq 0$. Тогда, в силу слабой полунепрерывность снизу функционала $T(u)$ на $\mathcal{D}^{1,2}(\Omega)$, вытекает 
\begin{align}
	&E_\lambda(u_\lambda) \leq \liminf_{m\to \infty} E_\lambda(u_m)=\hat{E_\lambda}, \label{ner3}\\
	&E'_\lambda(u_\lambda) \leq \liminf_{m\to \infty} E'_\lambda(u_m)=0,\label{ner4}
\end{align}
где $u_\lambda=r_0v_0$. Предположим, что в (\ref{ner4}) выполняется строгое неравенство. Тогда, учитывая $E'_\lambda(u_\lambda)<0$, найдется $t>1$ такое, что $E'_\lambda(tu_\lambda) =0$, $E''_\lambda(tu_\lambda)>0$ и $E_\lambda(tu_\lambda)<E_\lambda(u_\lambda)\leq \hat{E_\lambda}$. Получили противоречие. Следовательно в (\ref{ner4}) выполняется равенство. Тогда $u_\lambda \in \mathcal{N}_\lambda$, что  влечет  равенство в (\ref{ner3}). Таким образом, получили требуемое.

\hspace*{\fill}\rule{3mm}{3mm}\\

Завершим теперь доказательство утверждения (2). Пусть $\lambda>\lambda_{E,(p,q)}$. Тогда, т.к. $\lambda_{E,(p,q)}>\lambda_{(p,q)}$, то по лемме \ref{le1e} существует решение $u_\lambda$ проблемы Нехари (\ref{min1}). Поскольку $\lambda>\lambda_{E,(p,q)}$, то по лемме \ref{Rl} найдется  $u \in \mathcal{D}^{1,2}(D) \setminus 0$	такое, что $E_\lambda(u)<0$. Отсюда мы заключаем, что $\hat{E}_\lambda<0$ и, следовательно $E_\lambda(u_\lambda)<0$ и   $E''_\lambda(u_\lambda)\neq 0$. Таким образом, $u_\lambda$ удовлетворяет условию (\ref{N}) и следовательно уравнению (\ref{s1U}). Теорема доказана.

\hspace*{\fill}\rule{3mm}{3mm}\\

\section{Линейная неустойчивость}

В этом параграфе, мы приведем некоторые результаты об линейно неустойчивости стационарных  решений 
следующего уравнения: 
\begin{equation}
\label{p1}
\partial_t u=\Delta u + \lambda|u|^{p-2}u-|u|^{q-2}u,~~~(t,x) \in (0,\infty)\times D.
\end{equation}  
Здесь, как и выше, $D=\mathbb{R}^N$,  $D= \Omega$ или $D=\mathbb{R}^N\setminus \overline{\Omega}$. 
Мы рассматриваем решения с граничными условиями (\ref{bD}), (\ref{bDI}), (\ref{bI}) и с начальным условием 
\begin{equation}\label{CP}
	u|_{t=o}=v_0,
\end{equation}
Решение (\ref{p1})-(\ref{CP}) мы будем обозначать: $u(t;v_{0})$. Хорошо известно (см. напр. \cite{hernandez1, quitt}), что если $v_{0}\in 
\mathrm{L}^{\infty }(D)$,  то существует единственное классическое решение $u(\cdot ;v_{0})\in C^{2,1}([0,T )\times D )\cap C([0,T )\times \overline{D})$ задачи (\ref{p1})-(\ref{CP}) при некотором $T\equiv T(v_0)\in (0,+\infty)$.  Более того, если $v_0 \in C(\overline{D})$, то 
$u\in C((0 ,T),\mathcal{D}^{1,2}(D))$, и если к тому же $v_0 \in \mathcal{D}^{1,2}(D)$, то 
\begin{equation}
\int_{0}^{t}||u_{t}(s;v_0)||_{L^{2}}^{2}ds+E_{\lambda }(u(t;v_0))= E_{\lambda
}(v_0).  \label{ener1}
\end{equation}
В случае $T=+\infty$, решение $u(\cdot ;v_{0})$ называется глобальным.

Пусть $u_0 $ -- слабое, ограниченное в $D$ решение (\ref{s1U}). Рассмотрим соответствующую линеаризованную задачу  
\begin{equation}
\left\{ 
\begin{array}{ll}
&-L \psi:=-\Delta \psi -(\lambda (p-1) |u_0|^{p-2}-(q-2)|u_0|^{q-1})\psi =\mu \psi  , \\ 
\,&~~\psi \in W(D).%
\end{array}%
\right.  \label{line}
\end{equation}%
Тогда, при $2<\min\{p,q\}<+\infty$, не сложно показать, что существует минимальное собственное значение  $\mu_{1}$ задачи (\ref{line}) с неотрицательной собственной функцией  $\psi _{1}\in  W(D)$ (см. напр. \cite{giltrud, Reed}).

Решение  $u_0$ задачи (\ref{s1}) мы будем называть  
линейно неустойчивым стационарным решением  задачи (\ref{p1}),  если 
 минимальное собственное значение  $\mu_{1}$ оператора $-L$ отрицательно. Отметим, если рассмотреть 
возмущенное решение $v(t,x)=u_0(x)+w(t,x)$ задачи (\ref{p1}), то линеаризуя получаем   
$\partial_t w = L w$. Отсюда, при $w=e^{-t\mu_1}\psi_1$, мы убеждаемся, что в линейном приближении возмущение экспоненциально неустойчиво. 
 
\begin{lem}
\label{lemUn}  Пусть $\mathrm{N}\geq 1$, $p\neq q$, $p\geq 2, q\geq 2$, $\Omega $ -- ограниченная область, звездная  относительно начала координат  $\mathbb{R}^{N}$,  с $C^1$-гладкой границей $\partial \Omega$.  
 Тогда слабое, ограниченное в $D$ решение $u$ задачи (\ref{s1U}) является линейно неустойчивым стационарным решением  задачи (\ref{p1}), если 
\begin{description}
	\item[(i)] $D=\mathbb{R}^N$ и  $d^*(p,q)< 0$;
	\item[(ii)] $D=\Omega$,~ $\max\{2,q\}<p<2^*$ и $u$ неотрицательно в $D$;
	\item[(iii)] $D=\mathbb{R}^N\setminus \Omega$, $u$ неотрицательно в $D$, $2<p<q$, или $d^*(p,q)< 0$ и $0<q<p$.
\end{description}
\end{lem}

{\it Доказательство. }\, Заметим, что при условиях \textbf{(i)-(iii)}, по теоремам \ref{Th1}, \ref{Th12}, \ref{Th3}   слабое решение $u$ задачи (\ref{s1U}) удовлетворяет $E''(u)<0$. По минимаксному принципу Куранта-Фишера имеем  
\begin{equation}
\mu _{1}=\inf_{\psi \in \mathcal{D}^{1,2}(D) \setminus \{0\}}\frac{\displaystyle \int_{D
}\left (|\nabla \psi |^{2}-(\lambda (p-1) |u|^{p-2}-(q-2)|u|^{q-1})\psi ^{2}\right )\,dx}{\displaystyle\int_{D }|\psi |^{2}\,dx}
\label{Courant}
\end{equation}%
Положим $\psi =u$. Тогда 
\begin{equation}\label{Courant2}
\frac{\displaystyle\int_{D }\left (|\nabla u|^{2}-(\lambda (p-1) |u|^{p-2}-(q-1)|u|^{q-2})u^{2}\right )\,dx}{%
\displaystyle\int_{D }|u|^{2}\,dx}=\frac{E^{\prime \prime
}(u)}{\displaystyle\int_{D }|u|^{2}\,dx}.
\end{equation}%
Отсюда, поскольку $E^{\prime \prime }(u_{\lambda })<0$, то из \eqref{Courant} вытекает, что $\mu _{1}<0$. Лемма доказана. 
 
\hspace*{\fill}\rule{3mm}{3mm}\\

Отметим, что подход основанный на использовании функции расслоений $E(ru)$ позволяет дополнительно получить следующий результат об оценки сверху минимального собственного значения линеаризованной задачи  (\ref{line}). 

\begin{cor}
 Пусть $\mathrm{N}\geq 1$, $p\neq q$, $p\geq 2, q\geq 2$,  $D=\mathbb{R}^N$,  $D= \Omega$ или $D=\mathbb{R}^N\setminus \overline{\Omega}$, где $\Omega $ -- ограниченная область  с $C^1$-гладкой границей $\partial \Omega$. 	Пусть $u$ -- слабое, ограниченное в $D$ решение задачи (\ref{s1U}). Тогда для минимального 
собственного значения $\mu_1$ линеаризованной задачи  (\ref{line}) справедлива следующая оценка
$$
\mu_1\leq \frac{E^{\prime \prime
}(u)}{\displaystyle\int_{D }|u|^{2}\,dx}
$$
\end{cor}
{\it Доказательство}\, непосредственно следует из минимаксного принципа Куранта-Фишера \eqref{Courant}  и формулы \eqref{Courant2}.

\section{Об устойчивости основных состояний.}

Отметим, что теорема Деррика \cite{Derrick} дает более сильное утверждение, чем в лемме \ref{lemUn} при  условии \textbf{(i)}, а именно: при $D=\mathbb{R}^N$, $N\geq 3$ и $p>2, q>2$, если существует слабое, ограниченное в $\mathbb{R}^N$ решение $u$ задачи (\ref{s1U}), то оно является линейно неустойчивым стационарным решением   (\ref{p1}). Этот же результат, как это показано в \cite{Comech}, справедлив и при $N=1,2$, если $u \in H^\infty(\mathbb{R}^N)$. Однако, если $p,q \in (1,2]$, ситуация может отличается. Это связано, в частности и с тем, что в этом случае задача (\ref{s1U}) может обладать решениями с компактными носителями \cite{brezis1971, CortElgFelmer-2, Diaz-vol-1, Kaper1}.  Ниже мы покажем, что такие решения, в определенном смысле, являются устойчивыми, при $D=\mathbb{R}^N$.
Отличительной ситуация  возникает и, когда рассматривается задача в ограниченной области $D=\Omega$. В этом случае, существуют устойчивые решения задачи (\ref{s1U}). При $p,q \in (1,2]$ такой результат был получен в \cite{DIH}. В этом разделе, мы обобщаем этот результат, на случаи других показателей нелинейностей.

 В дальнейшим, решение $u_0$ задачи (\ref{s1U}) будем называть $\mathcal{D}^{1,2}(D)$-\textit{устойчивым} стационарным решением параболической задачи (\ref{p1})-(\ref{CP}), если $\forall \varepsilon >0
$, $\exists \delta >0$ такое, что
\begin{equation}
||u_{0}-u(t;v_0)||_{1}<\varepsilon, ~~%
\forall v_0 \in \mathcal{D}^{1,2}(D)\cap C(\overline{D}):~  ||u_0-v_0||_{1}<\delta,~~\forall t>0.
\end{equation}%
При фиксированном $\lambda>0$, множество всех основных состояний задачи (\ref{s1U}) мы будем обозначать $g_\lambda:=g_\lambda(D)$ и называть \textit{многообразие основных состояний}. 
\begin{prop}
Пусть $2<p< q<2^*$ и $\lambda>\lambda_{(p,q)}$, $D=\Omega$ -- ограниченная область в $\mathbb{R}^{N}$ с $C^1$-гладкой границей $\partial \Omega$. Тогда	$g_{\lambda}(D)$ -- ограниченное множество в $\mathcal{D}^{1,2}(D)$. 
\end{prop}
{\it Доказательство} этого утверждения  основывается на  рассуждении от противного и аналогично доказательству леммы \ref{le1e}.
   
Пусть $\delta >0$, обозначим 
$
V_{\delta }(g_\lambda):=\{v\in \mathcal{D}^{1,2}(D):\inf_{u\in 
g_\lambda}||u-v||_{1}<\delta \}\cap C(\overline{D}).
$
Многообразие основных состояний $g_\lambda$ мы будем называть  $\mathcal{D}^{1,2}(D)$-\textit{устойчивым} для параболической задачи (\ref{p1})-(\ref{CP}), если  $\forall \varepsilon >0
$ найдется  такое $\delta >0$, что
\begin{equation}
\inf_{u_{0}\in G_{\lambda }}||u_{0}-u(t;v_0)||_{1}<\varepsilon, ~~%
\forall v_0 \in V_{\delta }(G_{\lambda }),~~\forall t>0.
\end{equation}%

\begin{lem}\label{lem4}
 Пусть $\mathrm{N}\geq 1$, $1<q<p<2$ или $2<p<q<2^*$, $D=\Omega $,  где $\Omega $ -- ограниченная область в $\mathbb{R}^{N}$ с $C^1$-гладкой границей $\partial \Omega$. Тогда, если  $\lambda>\lambda_{(p,q)}$, то многообразие основных состояние $g_\lambda$ задачи (\ref{s1U}) является  $\mathcal{D}^{1,2}(D)$-\textit{устойчивым}. 
\end{lem}
{\it Доказательство}.  Рассмотрим многообразие основных состояний $g_{\lambda}$ задачи \eqref{s1U}. Заметим, что по теореме \ref{thmblO}, если $\lambda>\lambda_{(p,q)}$, то $E_{\lambda }''(u)>0$ для любого $u \in g_{\lambda}$. Поскольку $g_{\lambda}$ ограниченное множество в $\mathcal{D}^{1,2}(D)$, а отображения $E_{\lambda },E_{\lambda }'':\mathcal{D}^{1,2}(D)\rightarrow \mathbb{R}$ непрерывны, то найдется такое $\delta_0 >0$, что при всех $u\in
V_{\delta }(g_{\lambda })$ и $0<\delta <\delta _{0}$, будет выполняться $E_{\lambda }''(u)>0$.

Покажем, что $E_{\lambda }$ является функционалом Ляпунова в окрестности $V_{\delta }(g_{\lambda })$ при $0<\delta <\delta _{0}$.

\begin{prop}
\label{lemSt1}   Существует $\delta \in (0,\delta _{0})$ такое, что 
\begin{equation}
E_{\lambda }(u)>\hat{E_{\lambda }},~~\forall u\in
V_{\delta }(g_{\lambda })\setminus g_{\lambda }.
\end{equation}
\end{prop}

{\it Доказательство.\,} Предположим противное, т.е. при каждом
$\delta \in (0,\delta _{0})$ существует $u^{\delta }\in V_{\delta }(g_{\lambda })\setminus g_{\lambda }$ такое, что $E_{\lambda }(u^{\delta
})\leq\hat{E_{\lambda }}$. Тогда найдется последовательность  $u_n\in V_{\delta _{0}}(g_{\lambda })\setminus g_{\lambda }$ такая, что 
$$
\inf_{u \in g_\lambda}||u-u_n||_1 \to 0~~\mbox{при}~~ n\rightarrow
\infty 
$$
и
\begin{equation}\label{st2}
E_{\lambda }(u_n)\leq \hat{E_{\lambda }}\quad n=1,2,....
\end{equation}%
Из первого  вытекает, что существует последовательность  $v_n \in g_{\lambda }$ такая, что
\begin{equation}\label{strong}
	||v_n-u_n||_{1}\rightarrow
0~~\mbox{при}~~n\rightarrow \infty. 
\end{equation}
Заметим, что $(v_n)$ является минимизационной последовательностью \eqref{min1}, т.к. 
$E_\lambda(v_n)\equiv \hat{E}_\lambda$ при всех $n=1,2,...$. Следовательно мы можем применить те же рассуждения, что и использовались для минимизирующей последовательности при доказательстве леммы \ref{le1e}. Поэтому $(v_n)$ имеет предельную точку $u_\lambda \neq 0$ являющейся основным состоянием   задачи \eqref{s1U}. При этом, $v_m \to u_\lambda$ сильно в  $\mathcal{D}^{1,2}(\Omega)$ при $n\to \infty$. Это и  \eqref{strong} влечет, что $
u_n\rightarrow u_{\lambda }$ сильно в $\mathcal{D}^{1,2}(\Omega)$
при $n\rightarrow \infty $. 
Отметим, по построению $u_n$ не является основным состоянием \eqref{s1U}. Поэтому
\begin{equation*}
E_{\lambda }(r_{min}(u_n)u_n)>\hat{E_{\lambda }}\quad n=1,2,....
\end{equation*}%
Данное неравенство и \eqref{st2} влечет
\begin{equation}
1<r_{max}(u_n)<r_{\min }(u_n).  \label{st3}
\end{equation}%
Поскольку отображения $r_{max}(\cdot ),r_{\min }(\cdot ):\mathcal{D}^{1,2}(\Omega)\rightarrow 
\mathbb{R}$ непрерывны и $u_n\rightarrow u_{\lambda }$ в $\mathcal{D}^{1,2}(\Omega)$ при $%
n\rightarrow \infty $, то
\begin{equation*}
r_{\min }(u_n)\rightarrow r_{\min }(u_{\lambda })=1~~\mbox{при}~~
n\rightarrow \infty ,
\end{equation*}%
Здесь $r_{\min }(u_{\lambda })=1$, т.к. $E_{\lambda
}^{\prime \prime }(u_{\lambda })>0$.  Тогда из \eqref{st3} следует, что
\begin{equation*}
r_{max}(u_n)\rightarrow r_{\min }(u_{\lambda })=1~~\mbox{при}\quad n\rightarrow
\infty .
\end{equation*}%
Отсюда и, т.к. $E_{\lambda }''(r_{max}(u_n)u_n)\leq
0 $ и $E_{\lambda }''(r_{\min }(u_n)u_n)\geq 0$, мы получаем
\begin{equation*}
E_{\lambda }''(u_{\lambda })=0.
\end{equation*}%
Однако это противоречит тому, что $E_{\lambda }''(u)>0$ для всех $u \in g_\lambda$.

\hspace*{\fill}\rule{3mm}{3mm}

Для завершения доказательства леммы \ref{lem4} достаточно проверить следующее
\begin{prop}
\label{lemst2}  Для каждого $\varepsilon >0$ найдется $\delta \in (0,\delta _{0})$ такое, что   
\begin{equation}
\inf_{u \in g_\lambda}||u -u(t;v_0)||_{1}<\varepsilon,~ ~~\forall v_0\in
V_{\delta }(g_{\lambda }),~~\forall t>0.
\end{equation}
\end{prop}
{\it Доказательство.\,} Пусть  $%
\varepsilon \in (0,\delta _{0})$. Рассмотрим 
\begin{equation}\label{depsilon}
d_{\varepsilon }:=\inf \{E_{\lambda }(w):w\in \mathcal{D}^{1,2}(\Omega),~\inf_{u \in g_\lambda}||u-w||_{1}=\varepsilon \}.
\end{equation}%
Тогда $d_{\varepsilon }>\hat{E_{\lambda }}$. Действительно, предположим противное, т.е. существует последовательность 
$w^{n}\in \mathcal{D}^{1,2}(\Omega)$ такая, что $\inf_{u \in g_\lambda}||u-w^{n}||_{1}=\varepsilon $ и $E_{\lambda }(w^{n})\rightarrow \hat{%
E_{\lambda }}$. Поскольку множество $g_\lambda$ ограниченно в $\mathcal{D}^{1,2}(\Omega)$, то  $(w^{n})$ также ограниченно в $\mathcal{D}^{1,2}(\Omega)$. Тогда по теореме вложения Соболева существует $v_0\in \mathcal{D}^{1,2}(\Omega)$ и подпоследовательность (снова обозначаемая 
$(w^{n})$) такая, что $w^{n}\rightarrow v_0$ слабо в $\mathcal{D}^{1,2}(\Omega)$ и сильно в $L_{\gamma},~ 1<\gamma<2^*$. Из слабой полунепрерывности снизу функционала $||u||_{1}^{2}$ на $\mathcal{D}^{1,2}(\Omega)$, вытекает $\hat{E_{\lambda }}\geq
E_{\lambda }(v_0)$ и $\inf_{u \in g_\lambda}||u-v_0||_{1}\leq \varepsilon $. По
предложению \ref{lemSt1} это возможно только, если $v_0 \in g_\lambda$. Но тогда равенство $\hat{E}_{\lambda }%
=E_{\lambda }(v_0)$ влечет, что $w^{n}\rightarrow v_0$ сильно в $%
\mathcal{D}^{1,2}(\Omega)$. Отсюда   $\varepsilon =\inf_{u \in g_\lambda}||u-w^{n}||_{1}\rightarrow \inf_{u \in g_\lambda}||u-v_0||_{1}$. Но тогда $v_0\notin  g_{\lambda }$. Получили противоречие.

Поскольку $\hat{E_{\lambda }}<d_{\varepsilon }$, то $\hat{E_{\lambda }}<d_{\varepsilon }-\sigma$ для некоторого $\sigma >0$. Из непрерывности отображения $E_{\lambda }:\mathcal{D}^{1,2}(\Omega)\rightarrow 
\mathbb{R}$ вытекает, что найдется такое $%
\delta \in (0,\varepsilon )$, что
\begin{equation}
E_{\lambda }(w)<d_{\varepsilon }-\sigma, \quad\forall w\in V_{\delta
}(g_{\lambda })\subset V_{\varepsilon }(g_{\lambda }).  \label{eqst3}
\end{equation}%
Для доказательства предложения, остается убедится, что для любого $v_0\in V_{\delta }(g_{\lambda })$ решение $
u(t,v_0)$ остается в $V_{\varepsilon }(g_{\lambda })$ при всех $t>0$.
Предположим противное. Тогда, учитывая, что $u(t,v_0)\in
C((0,T),\mathcal{D}^{1,2}(\Omega))$, найдется $t_{0}>0$ такое, что $\inf_{u \in g_\lambda}||u-u(t_{0},v_0)||_{1}=\varepsilon $. Но это, в силу \eqref{depsilon}, влечет
\begin{equation*}
d_{\varepsilon }\leq E_{\lambda }(v(t_{0},v_0)).
\end{equation*}%
С другой стороны, из (\ref{ener1}) вытекает, что $E_{\lambda
}(v(t_{0},v_0))\leq E_{\lambda }(v_0)$. Таким образом, учитывая \eqref{eqst3}, имеем
\begin{equation*}
d_{\varepsilon }\leq E_{\lambda }(v(t_{0},v_0))\leq E_{\lambda
}(v_0)<d_{\varepsilon }-\sigma.
\end{equation*}%
Получили противоречие. Предложение  доказано.

\hspace*{\fill}\rule{3mm}{3mm}\\

Пусть $D=\Omega$. Если неотрицательное решение $u$ задачи (\ref{s1U}) удовлетворяет  
\begin{equation}
\frac{\partial u}{\partial \nu }=0~~\mbox{на}~~\partial \Omega,
\label{NV}
\end{equation}%
то говорят, что $u$ является \textit{решением с компактным носителем}. Отметим, что по стандартной теории регулярности решений эллиптических уравнений, если $u_{\lambda^{\ast}} \in W(D)$ слабое решение 
(\ref{s1U}), то $u \in C^2(D)\cap C^{1,\kappa}(\overline{D})$ для  $\kappa \in (0,1)$. 

В \cite{ ilDr, IlEg} доказано следующее утверждение
\begin{lem}\label{lem5}
\label{Th2} Пусть  $\mathrm{N}\geq 3$, $1<q<p<2$ и $d^*(p,q)>0$, $D= \Omega$  -- ограниченная область, звездная относительно начала координат  $\mathbb{R}^{N}$  с $C^1$-гладкой границей $\partial \Omega$.
Тогда найдется такое $\lambda^{\ast }\equiv\lambda^{\ast }_{(p,q)}>\lambda_{(p,q)}$, что
при всех $\lambda \geq \lambda^{\ast }$ задача (\ref{s1U})  обладает решением $u_\lambda^c$ с компактным носителем. Более того,  при $\lambda = \lambda^{\ast }_{(p,q)}$, решение с компактным носителем $u_{\lambda^{\ast}}^c$ является основным состоянием (\ref{s1U}). При этом, $u_\lambda^c \in C^2(D)\cap C^{1,\kappa}(\overline{D})$ для  $\kappa \in (0,1)$ и  $u_\lambda^c\geq 0$ в $D$.  
\end{lem}
Отметим, что $\{(p,q) \in \mathbb{R}^2: 1<q<p<2\} \cap \{(p,q) \in \mathbb{R}^2: d^*(p,q)>0\}=\emptyset$
при $N=1,2$.

Нам понадобится следующий результат из \cite{serrin}

\begin{lem}
\label{lem:3} (Serrin-Zou) Пусть $\mathrm{N}\geq 2$, $1<q<p<2$, $D= \mathbb{R}^N$. Тогда любое решение  $u$ задачи (\ref{s1U}) имеет компактный носитель. Более того, для каждой отдельной связанной компоненты $\Xi $ в открытом носителе $\Theta :=\{x\in \mathbb{R}^{N}:u(x)>0\} $ выполняется следующее: 
\begin{enumerate}
\item $\Xi$ является шаром;
\item $u$ -- радиально симметричная функция относительно центра шара $\Xi $.
\end{enumerate}
\end{lem}

Пусть $\mathrm{N}\geq 3$, $1<q<p<2$. Рассмотрим задачу (\ref{s1U}) при $D= \mathbb{R}^N$ и $\lambda=1$. Тогда по лемме \ref{ThmEx1} существует неотрицательное классическое решение $u^c$ этой задачи, являющиеся основным состоянием (\ref{s1U}). Из леммы \ref{lem:3}, поскольку $u^c$ -- основное состояние, несложно заключить, что носитель supp($u^c$) состоит из одной компоненты,  являющейся шаром, без ограничения общности с центром в нуле, с некоторым радиусом $R_{(p,q)}>0$, т.е. supp($u^c)=B_{R_{(p,q)}}$. При этом, $u^c$ -- радиально-симметричная функция. Отметим, что $u^c$ является классическим решением   задачи (\ref{s1U}) при $D= B_{R_{(p,q)}}$. Отсюда, как это показано в  \cite{Kaper1,Kaper2}, вытекает, что  $u^c$ -- единственное положительное решение этой задачи, а радиус $R_{(p,q)}$ задается однозначным образом. С другой стороны, по лемме \ref{lem5}, если дополнительно выполняется $d^*(p,q)>0$, то существует $\lambda^{\ast }_{(p,q)}(B_{R_{(p,q)}})$ такое, что  (\ref{s1U}) обладает основным  состоянием  $u_{\lambda^{\ast}}^c$ с компактным носителем. Используя отмеченную единственность   $u^c$ и то, что  $u^c$ -- основное состояние задачи (\ref{s1U}) при $D= \mathbb{R}^N$, несложно показать, что $\lambda^{\ast }_{(p,q)}(B_{R_{(p,q)}})=1$ и $u^c=u_{\lambda^{\ast}}^c$. Отсюда, в частности, получаем  

\begin{cor}\label{clem4}
 Пусть $\mathrm{N}\geq 3$, $D=B_{R_{(p,q)}} $, $1<q<p<2$ и $d^*(p,q)>0$. Тогда многообразие основных состояний $g_{\lambda^{\ast}}(B_{R_{(p,q)}})$ задачи (\ref{s1U}) с $\lambda=1$ состоит из единственного решения
$u^c$, при этом $u^c$ имеет компактный носитель.
\end{cor}
Любая функция  $w$ из $\mathcal{D}^{1,2}(\Omega)$
может быть продолжена в $\mathbb{R}^{N}$ по формуле: 
\begin{equation}
\left\{ 
\begin{array}{ll}
\tilde{w}=w &\mbox{в}~\Omega , \\ 
\tilde{w}=0 & \mbox{в}~\mathbb{R}^{N}\setminus \Omega ,%
\end{array}%
\right.  \label{expan}
\end{equation}%
Тогда $\tilde{w}\in \mathcal{D}^{1,2}(\mathbb{R}^{N})$ и, в этом смысле 
$\mathcal{D}^{1,2}(\Omega)\subset \mathcal{D}^{1,2}(\mathbb{R}^{N})$.

Из леммы \ref{lem4} и следствия \ref{clem4}, имеем
\begin{lem}\label{lemF}
Пусть $\mathrm{N}\geq 3$, $D= \mathbb{R}^N$, $1<q<p<2$ и $d^*(p,q)>0$. Тогда решение с компактным носителем $u^c$ задачи  
(\ref{s1})-(\ref{bDI}) является устойчивым стационарным состоянием параболической задачи (\ref{p1})-(\ref{CP}) в следующем смысле:  
$\forall \varepsilon >0
$, $\exists \delta >0$ такое, что
\begin{equation}
||u^c-u(t;v_0)||_{1}<\varepsilon, ~~%
\forall v_0 \in \mathcal{D}^{1,2}(B_{R_{(p,q)}})\cap C(\overline{D}):~  ||u^c-v_0||_{1}<\delta,~~\forall t>0.
\end{equation}%
\end{lem}
Из трансляционной инвариантности уравнений  (\ref{s1}), (\ref{p1}), очевидно $u^c(\cdot +y)$ при любом $y \in \mathbb{R}^N$ также является устойчивым стационарным состоянием параболической задачи (\ref{p1})-(\ref{CP}) в аналогичном смысле, что и в лемме \ref{lemF}. 
В частности, рассмотрим 
\begin{equation*}
	\mathcal{M}_\delta:=\{v \in \mathcal{D}^{1,2}(\mathbb{R}^N):\exists y \in \mathbb{R}^N,~v(\cdot +y) \in  \mathcal{D}^{1,2}(B_{R_{(p,q)}}),~\sup_{y \in  \mathcal{O}_v} ||u^c(\cdot)-v(\cdot +y)||_1<\delta \},
\end{equation*}
где $\delta>0$, $\mathcal{O}_v:=\{y \in \mathbb{R}^N:~v(\cdot +y) \in  \mathcal{D}^{1,2}(B_{R_{(p,q)}})\}$ для $v \in \mathcal{D}^{1,2}(\mathbb{R}^N)$. Тогда справедлив следующий результат о существовании глобальных решений
\begin{cor}\label{clem5}
 Пусть $\mathrm{N}\geq 3$, $D= \mathbb{R}^N$, $1<q<p<2$ и $d^*(p,q)>0$. Существует $\delta_0>0$ такое, что при всех $v_0 \in \mathcal{M}_{\delta_0}$, решение $u(t;v_0)$ задачи (\ref{p1})-(\ref{CP})  
является глобальным и ограниченным при всех $t>0$.
\end{cor}

\section{Заключительные замечания}

Несложно заметить, что теорема \ref{thmblO} о существовании решений, в случае ограниченной области $D=\Omega$, обобщается на уравнения (\ref{s1}) с ненулевой массой, по крайне мере при малых $m>0$. 
Учитывая теорему \ref{Th12}, открытым остается вопрос о существовании решений (\ref{s1}) при $D=\Omega$,
когда $2^*\leq p<q$, $0<p<q\leq 1$,  $0<q<p\leq 1$.

Нам не известны общие результаты о существовании решений (\ref{s1U}), если $D=\mathbb{R}^N\setminus \Omega$, при наличия двух показателей нелинейности $p,q$. Представляется, что в целом данная проблема остается открытой.  Основные результаты по этой теме имеются, насколько нам известно, только для задач с одним показателем нелинейности (см напр.  \cite{BahriLions, davila, pucci} и приведенные там ссылки).   

Наличие двух и более стационарных точек функции расслоений, как для (\ref{ffib}) при $1<q<p<2$ или $2<p<q$, обычно влечет кратность  решений (multiplicity solutions) и существованию ветвей решений с    бифуркационной точкой типа поворота (см. \cite{Amb, DrP, ilyas, ilyas2}). Однако утверждение $(2^o)$	 теоремы \ref{Th1}, дает основание выдвинуть гипотезу о том, что такое поведение для  (\ref{s1U}) при $D=\mathbb{R}^N$ отсутствует. Аналогичные гипотезы возникают и на основание утверждений $(2^o)$	 теорем \ref{Th12}, \ref{Th3} для задач (\ref{s1U}) при $D=\Omega$ и $D=\mathbb{R}^N\setminus \Omega$. 


Отметим, что Деррик в \cite{Derrick}, формулируя результат  об отсутствие устойчивости для локализованных решений, сделал  на этом основание заключение о том, что данное свойство является препятствием  для интерпретации рассматриваемых решений как частиц. Он предложил в \cite{Derrick} несколько способов модификации рассматриваемых моделей с целью получения устойчивых локализованных решений. Из леммы \ref{lemF} и следствия \ref{clem5} можно сделать предположение, что уравнения с нелипшицевыми нелинейностей могут рассматриваться как один из способов получения моделей с устойчивыми локализованными решениями.

\end{document}